\documentclass[12pt]{article}
\usepackage{graphicx, natbib}

\newcommand{\tr}{^{\prime}}

\def\b#1{\mbox{\boldmath $#1$}}    % - \b: grassetto in formula
\def\bl#1{\mbox{\footnotesize \boldmath {$#1$}}} % - \bl: grassetto per apice
\renewcommand{\th}{\theta}
\newcommand{\pa}{\partial}
\newcommand{\al}{\alpha}
\newcommand{\be}{\beta}
\newcommand{\de}{\delta}
\newcommand{\la}{\lambda}
\newcommand{\La}{\Lambda}
\newcommand{\ga}{\gamma}
\newcommand{\si}{\sigma}

\newcommand{\diag}{{\rm diag}}    % - \diag: simbolo diag
\usepackage{color}
\definecolor{colorS}{rgb}{0.7695312,0.0468750,0.1601562}
\begin{document}

\title{\vspace*{-1.5cm} Nested hidden Markov chains for modeling\\
dynamic unobserved heterogeneity\\ in multilevel longitudinal data}
\author{Francesco Bartolucci and Monia Lupparelli}
%% Use \authorrunning{Short Title} for an abbreviated version of
%% your contribution title if the original one is too long
%\institute{Francesco Bartolucci \at University of Perugia
%\email{bart@stat.unipg.it} \and Monia Lupparelli \at University of
%Bologna \email{monia.lupparelli@unibo.it}
%}
%
% Use the package "url.sty" to avoid
% problems with special characters
% used in your e-mail or web address
%
\maketitle

\begin{abstract}
In the context of multilevel longitudinal data, where sample units are collected in clusters, an important aspect that should be accounted for is the unobserved heterogeneity between sample units and between clusters. For this aim we propose an approach based on nested hidden (latent) Markov chains, which are associated to every sample unit and to every cluster. The approach allows us to account for the mentioned forms of unobserved heterogeneity in a dynamic fashion; it also allows us to account for the correlation which may arise between the responses provided by the units belonging to the same cluster. Given the complexity in computing the manifest distribution of these response variables, we make inference on the proposed model through a composite likelihood function based on all the possible pairs of subjects within every cluster. The proposed approach is illustrated through an application to a dataset concerning a sample of Italian workers in which a binary response variable for the worker receiving an illness benefit was repeatedly observed.\vspace*{0.5cm}

\end{abstract}\newpage
\section{Introduction}\label{se.intro}
In modeling longitudinal data, it is common to account for the unobserved heterogeneity between sample units, that is, the heterogeneity that cannot be explained on the basis of the observable covariates \citep{digg:02,hsia:03,free:04,fitz:09}. This is normally accomplished by the introduction of latent variables or random effects. For instance, a typical approach consists of associating a random intercept to every sample unit which affects the distribution of each occasion-specific response in the same fashion. This allows us to account for a form of {\em time-constant} unobserved heterogeneity which is due to unobservable covariates and related factors.

More recent approaches for longitudinal data are based on allowing for a form of {\em time-varying} unobserved heterogeneity, relaxing in this way the assumption that the effect of unobservable covariates on the response variables is constant in time. This is sensible in many applied contexts, especially in the presence of long panels and with a limited set of observable covariates. Among these time-varying approaches, it is worth mentioning the one described in \cite{hei:08}, which is based on random effects having an AR(1) structure, and that proposed by \cite{bart:farc:09}, which is based on a hidden (latent) Markov chains for capturing the unobserved heterogeneity in a dynamic fashion. For a comparison between the two approaches see \cite{bart:bacci:penn:2010}.

The above considerations are obviously pertinent when we deal with multilevel longitudinal data, where sample units are collected in clusters, with the addiction that it is also appropriate modeling the unobserved heterogeneity between clusters and the correlation between the responses provided by the units in the same cluster. Note that multilevel longitudinal data are more and more easily encountered in socio-economic contexts. In particular, the dataset motivating this paper, that will be described in detail in the following, concerns a sample of workers (sample units) in different firms (clusters), who are longitudinally observed. As response we have a binary variable equal to 1 if the employee receives illness benefits in a certain year and to 0 otherwise. Datasets having a similar structure are nowadays available, for instance, in educational contexts, where students are collected in classes and are followed for a certain number of years of schooling. In these datasets we typically have a limited set of observable covariates and the need arises for an appropriate modeling of the unobserved heterogeneity between both sample units and clusters.

For the aim described above, we propose an approach based on nested hidden Markov chains which may be seen as an extension of the approach proposed by \cite{bart:farc:09} for longitudinal data. In particular, we associate a first-order homogeneous hidden Markov chain to every sample unit and to every cluster. The time-specific realizations of these two chains go to affect the distribution of the response variables together with the covariates observed at unit and cluster levels. Coming back to the above example about the sample of employees, the different states of the unit-level Markov chain correspond to different levels of the residual (not explained by the unit-level observable covariates) tendency to require an illness benefit by an employee. A similar interpretation may be found for the different states of the cluster-level Markov chain, which affect the behavior of the employees in the same firm. Moreover, the possibility that the unit-level state changes may be due to events of the employee's life that are not recorded in the dataset, such as a sudden worsening of his/her health status. Similarly, a change in the cluster-level state may be due to events about the firm, such as the change of the management. In any case, we can test if the latent effects are indeed dynamic or not on the basis of the dataset at hand.

The proposed approach may be cast in the literature about latent Markov (LM) models for longitudinal data, as described by \cite{bar:far:pen:10}. It is worth noting that other multilevel extensions of the latent (or hidden) Markov approach for longitudinal data are available in the literature. We mention, in particular, the extensions proposed by \cite{bart:lupp:mont:09} and \cite{bart:penn:vitt:11}. About multilevel extensions see also \cite{aspa:muth:08} and about related models including random effects, but not in a context of  analysis of multilevel data, see \cite{vand:lang:90}, \cite{Altman:07}, and \cite{maru:11}. In these cases the effects (fixed or random) associated to every cluster are time-constant. However, an extension in which these effects are time-varying has not been proposed yet, at least to our knowledge.

Under the proposed model, the manifest distribution of the response variables is computationally intractable in most applications. Therefore, to make inference on the model we exploit an approach based on a composite likelihood \citep{lindsay:88,cox:reid:04}, which is computed on the basis of the joint distribution of the response variables for each pair of subjects in the same cluster. A similar approach was followed by \cite{ren:mol:gey:04} to deal with a multilevel probit model; for applications of this inferential approach to similar contexts, see \cite{hjo:var:08} and \cite{var:cza:10}. In particular, we show how to compute the pairwise likelihood by using the same recursion exploited by \cite{baum:et:al:70} to deal with hidden Markov models and how to maximize this likelihood by an Expectation-Maximization (EM) algorithm similar to the one they suggest and implemented along the same lines as in \cite{bart:farc:09}. We also show how to obtain standard errors for the parameter estimates and how to make model selection on the basis of the composite likelihood information criterion (CLIC) developed by \cite{var:vid:05}. An \verb"R" implementation of the functions used for the estimation of the model in the presence of binary response variables is available to the reader upon request.

The paper is organized as follows. In the next section we briefly review the LM model with covariates \citep{bart:farc:09} and its maximum likelihood estimation. Section 3 illustrates the proposed multilevel extension dealing with the case of continuous and binary response variables. Pairwise likelihood inference for this model is described in Section 4. In Section 5 we illustrate the model by an application based on the dataset concerning the sample of workers mentioned above. Finally, in Section 6 we draw the main conclusions.
\section{Using hidden Markov chains for modeling unobserved heterogeneity}\label{sec.LM}
Consider a panel of $n$ subjects observed at $T$ occasions and let
$Y_i^{(t)}$ denote the response variable of interest for subject $i$
at occasion $t$, $i=1,\ldots,n$, $t=1,\ldots,T$, and let $\b Z_i^{(t)}$ be the corresponding column vector of covariates, which may also include the lagged responses. In the context of our application, the response variables are binary, although the LM model may be also applied to variables having a different nature.

In the following, we outline how to model these data accounting for unobserved heterogeneity in a dynamic fashion, by introducing a hidden Markov chain, as suggested by \cite{bart:farc:09}.
\subsection{Model assumptions}
We assume that, for $i=1,\ldots,n$, the response variables $Y_i^{(1)},\ldots,Y_i^{(T)}$ are conditionally independent given the covariate vectors $\b Z_i^{(1)},\ldots,\b Z_i^{(T)}$ and a latent process ${\b V}_i=(V_{i}^{(1)},\ldots,V_{i}^{(T)})$, which follows a first-order homogeneous Markov chain and is independent of the covariates.

This chain has $k$ states, labeled from 1 to $k$, with initial and transition probabilities
\begin{eqnarray*}
\pi_v&=&p(V_i^{(1)}=v),\quad v=1,\ldots,k,\\
\pi_{v|\bar{v}}&=&p(V_i^{(t)}=v|V_i^{(t-1)}=\bar{v}),\quad
t=2,\ldots,T, \: \bar{v},v=1,\ldots,k.
\end{eqnarray*}
Note that, in the above definitions, $v$ refers to the current state, whereas $\bar{v}$ refers to the previous one. This convention will be used throughout the paper. Moreover, the initial probabilities are collected in the $k$-dimensional column vectors $\b\pi$,
whereas the transition probabilities are collected in the $k\times k$ transition matrix $\b\Pi$. Note that these probabilities are the same for all sample units and, in particular, the transition probabilities are time homogenous. Moreover, in order to make the model more parsimonious, different constraints may be imposed on the matrix $\b\Pi$; see also \cite{bart:06}. For instance, we may assume that this matrix is tridiagonal, with constant off-diagonal elements, so that with $k=3$ we have
\begin{equation}
\b\Pi = \pmatrix{
1-\rho & \rho & 0\cr
\rho & 1-2\rho & \rho\cr
0 & \rho & 1-\rho
}\label{eq:tri},
\end{equation}
where $\rho$ is a parameter between 0 and 0.5 to be estimated.

For subject $i$ at occasion $t$, the latent variable $V_i^{(t)}$ corresponds to the level of the unobservable characteristic of interest. The way in which this characteristic affects the corresponding response variable $Y_i^{(t)}$ depends on the assumed {\em measurement model}. For instance, in the case of continuous response variables, it is natural to formulate the following assumption on the conditional distribution of $Y_i^{(t)}$ given $V_i^{(t)}$ and $\b Z_i^{(t)}$:
\[
Y_i^{(t)}|V_i^{(t)}=v,\b Z_i^{(t)}=\b z\sim N(\be_v+\b z\tr\b\de,\si^2),
\]
where $\be_v$ is an intercept related to the latent state and $\b\de$ is a vector of regression coefficients. Obviously, these parameters, including the variance $\si^2$, can be estimated together with the above initial and transition probabilities.

With binary response variables, instead, it is natural to assume that
\[
Y_i^{(t)}|V_i^{(t)}=v,\b Z_i^{(t)}=\b z\sim Bern(\psi_i^{(t)}(v,\b z)),
\]
where
\[
\log\frac{\psi_i^{(t)}(v,\b z)}{1-\psi_i^{(t)}(v,\b z)}=\be_v+\b z\tr\b\de,
\]
with $\psi_{hi}^{(t)}(v,\b z)$ corresponding to the conditional ``probability of success'', that is $\psi_{hi}^{(t)}(v,\b z)=p(Y_i^{(t)}=1|V_i^{(t)}=v,\b Z_i^{(t)}=\b z)$.

The above approach may be extended to response variables having a different nature, even ordinal variables, and also to multivariate contexts, in which we observe more response variables at each time occasions. We refer the reader to \cite{bart:farc:09} for details on the resulting LM model.
\subsection{Maximum likelihood estimation}\label{sec:EM}
When we deal with an observed sample, for $i=1,\ldots,n$ we have an observed response configuration ${\b y}_i = (y_i^{(1)},\ldots, y_i^{(T)})$ and an observed sequence of covariates vectors $\b z_1^{(1)},\ldots,\b z_i^{(T)}$; we collect these covariates in the unique vector $\b z_i$ (for all time occasion). In order to perform maximum likelihood estimation of the above model on the basis of these data, the need arises of computing the {\em manifest distribution} of $\b y_i$ given $\b z_i$, that is,
\begin{equation}
p(\b y_i|\b z_i)=
\sum_{\bl v}p(\b y_i|\b V_i=\b v,\b z_i)p(\b V_i=\b v),\label{eq:mani}
\end{equation}
where the sum $\sum_{\bl v}$ is over all the possible
configurations $\b v = (v_i^{(1)},\ldots, v_i^{(T)})$ of the latent process $\b V_i$.

Efficient computation of the probability in (\ref{eq:mani}) may be performed by exploiting a forward recursion available in the hidden Markov literature \citep[see][]{baum:et:al:70,levi:rabi:sond:83,macd:zucc:97}. As in \citet{bart:06} and \cite{bart:farc:09}, it is convenient to express this recursion by using the matrix notation on the basis of the initial probability vectors $\b\pi$ and transition matrix $\b\Pi$. For this aim, consider the column vector $\b q_i^{(t)}$ with elements
\[
p(y_i^{(1)},\ldots,y_i^{(t)},V_i^{(t)} = v,\b z_i^{(1)},\ldots,\b z_i^{(t)}),
\quad v=1,\ldots,k.
\]
This vector may be recursively computed as follows:
\begin{equation}\label{eq:recursion}
\b q_i^{(t)} = \left\{
\begin{array}{lll}
\diag({\b m}_i^{(1)})\b\pi, && \hbox{if } t = 1,  \\
\diag({\b m}_i^{(t)})\b\Pi\tr\b q_i^{(t-1)}, && \hbox{otherwise},
\end{array}%
\right.
\end{equation}
where $\b m_i^{(t)}$ is the column vector with elements $p(y_i^{(t)}|V_i^{(t)}=v,\b z_i^{(t)})$, for $v =1,\ldots,k$, which is defined on the basis of the assumed measurement model. Once this recursion has been performed for $t = 1,\ldots,T$, we may obtain $p(\b y_i)$ as the sum of the elements of the vector $\b q_i^{(T)}$.

Maximum likelihood estimation is performed by maximizing the log-likelihood $\ell(\b\th) = \sum_i\log[p(\b y_i|\b z_i)]$, where $\b\th$ denotes the vector of all model parameters. We maximize this function by an EM algorithm \citep{baum:et:al:70,demp:lair:rubi:77}, which is based on the {\it complete data log-likelihood} denoted by $\ell^*(\b \th)$, that is, the log-likelihood that we could compute if we knew the latent state of each subject at every occasion.

The EM algorithm alternates two steps (E and M) until convergence:
the \emph{E}-step computes the conditional expectation of
$\ell^*(\b\th)$, given the observed data and the current value of
$\b\th$, using recursions similar to the one illustrated above; the
\emph{M}-step maximizes this expected value with respect to
$\b\th$, so that this parameter vector results updated. The latter may require simple iterative algorithms of Newton-Raphson type. A detailed description of this EM algorithm is available in \cite{bart:farc:09}.
\section{Proposed multilevel extension}\label{sec.ext}
In the context of multilevel longitudinal data, the $n$ sample units are grouped, according to some criteria, in $H$ clusters of size $n_1,\dots,n_H$. Then, for each subject $i$ in cluster $h$, data are available at $T$ consecutive occasions. In particular, we denote by $Y_{hi}^{(t)}$ the corresponding  response variable and by $\b Z_{hi}^{(t)}$ the corresponding column vector of covariates, where $h=1,\ldots,H$, $i=1,\ldots,n_h$, and $t=1,\ldots,T$. Moreover, by $\b X_h^{(t)}$, with $h=1,\ldots,H$ and $t=1,\ldots,T$, we denote column vectors of cluster-level covariates, which may be time-varying.

In the following we show how multilevel longitudinal data, having the structure described above, may be analyzed by an extension of the approach outlined in Section \ref{sec.LM}.
\subsection{Model assumptions}
Our extension assumes the existence of a latent process $\b U_h=(U_{h}^{(1)},\ldots,U_{h}^{(T)})$ for each cluster $h$, $h=1,\ldots,H$, and a latent process $\b V_{hi}=(V_{hi}^{(1)},\ldots,V_{hi}^{(T)})$ for each subject $i$, $i=1,\ldots,n_h$, in the cluster. Both processes follow a first-order homogeneous Markov chain with $k_1$ states at cluster level and $k_2$ at individual level. These processes are assumed to be independent each other and also independent of the unit- and cluster-level covariates. Moreover, extending the assumptions  formulated in Section \ref{sec.LM}, we impose that, for every sample unit $hi$ (unit $i$ in cluster $h$), the response variables $Y_{hi}^{(t)}$ are conditionally independent given $\b U_h$, $\b V_{hi}$ and the corresponding covariates. This implies that the response vectors for two subjects in the same cluster are conditionally independent given $\b U_h$, but they are not marginally independent. This marginal independence holds for subjects belonging to two different clusters.

The initial and the transition probabilities of each cluster-level latent process are denoted by
\begin{eqnarray*}
\la_u&=&p(U_{h}^{(1)}=u),\quad u=1,\ldots,k_1\\
\la_{u|\bar{u}}&=&p(U_{h}^{(t)}=u|U_{h}^{(t-1)}=\bar{u}),\quad
t=2,\ldots,T,\:\bar{u},u=1,\ldots,k_1,
\end{eqnarray*}
and are collected in
the vector $\b\la$ and in the transition matrix $\b\La$. Moreover, for the unit-level latent processes $\b U_{hi}$, we substantially adopt the same notation as in Section \ref{sec.LM}, and then we let $\pi_v=p(V_{hi}^{(1)}=v)$ and $\pi_{v|\bar{v}}=p(V_{hi}^{(t)}=v|V_{hi}^{(t)}=\bar{v})$; these initial and transition probabilities are still collected in the vector $\b\pi$ and in the matrix $\b\Pi$, respectively.

Finally, about the conditional response probabilities, the same considerations expressed in Section \ref{sec.LM} still holds. Then, in the case of continuous response variables we may assume that:
\[
Y_{hi}^{(t)}|U_h^{(t)}=u,V_{hi}^{(t)}=v,\b X_h^{(t)}=\b x,\b Z_{hi}^{(t)}=\b z
\sim N(\al_u+\be_v+\b x\tr\b\ga+\b z\tr\b\de,\si^2),
\]
where $\al_u$ is an intercept related to the cluster-level latent state, $\be_v$ is an intercepts related to the unit-level latent state, and $\b\ga$ and $\b\de$ are corresponding vectors of regression coefficients.

With binary response variables, instead, it is natural to assume that
\begin{equation}
Y_{hi}^{(t)}|U_h^{(t)}=u,V_{hi}^{(t)}=v,\b X_h^{(t)}=\b x,\b Z_{hi}^{(t)}=\b z
\sim Bern(\psi_{hi}^{(t)}(u,v,\b x,\b z)),
\label{eq:cov_psi}
\end{equation}
where
\[
\log\frac{\psi_{hi}^{(t)}(u,v,\b x,\b z)}{1-\psi_{hi}^{(t)}(u,v,\b x,\b z)}=
\al_u+\be_v+\b x\tr\b\ga+\b z\tr\b\de,
\]
with parameters having the same interpretation as above.
\subsection{Manifest distribution}
When we observe a set of multilevel longitudinal data, we have a sequence of response $\b y_{hi}=(y_{hi}^{(1)},\ldots,y_{hi}^{(T)})$ for every sample unit $hi$, with $h=1,\ldots,H$, $i=1,\ldots,n_h$. We denote by $\b y_h$ the vector obtained by collecting the responses of all subjects in cluster $h$, that is $y_{hi}^{(t)}$ for $i=1,\ldots,n_h$ and $t=1,\ldots,T$. Similarly, we observe the vectors of unit-level covariates $\b z_{hi}^{(t)},\ldots,\b z_{hi}^{(T)}$; these covariates are collected in the unique vector $\b z_{hi}$ when referred to the unit $hi$ (for all time occasions) and in the vector $\b z_h$ when referred to all units in the same cluster $h$. Finally, for every cluster $h$, we observe the vectors of cluster-level covariates $\b x_h^{(t)}$, which are collected in the unique vector $\b x_h$ (for all time occasions).

Under the above assumptions, the manifest probability of $\b y_h$ given $\b x_h$ and $\b z_h$ has the following expression:
\begin{eqnarray*}
p(\b y_h|\b x_h,\b z_h)&=&\sum_{\bl u} p(\b U_h=\b u)\\
&&\times\prod_{i=1}^{n_h}\bigg[\sum_{\bl v} p(\b y_{hi}|\b U_h=\b u,\b V_{hi}=\b v,\b x_h,\b z_h)p(\b V_{hi}=\b v)\bigg],
\end{eqnarray*}
where the sum $\sum_{\bl u}$ is over all the possible configurations of the latent process $\b U_h$ and $\sum_{\bl v}$ is over all the possible configurations of $\b V_{hi}$.

For the cases in which computing $p(\b y_h|\b x_h,\b z_h)$ is feasible, estimation of the model parameters can be performed by maximizing the log-likelihood $\ell(\b\th) = \sum_h \log[p(\b y_h|\b x_h,\b z_h)]$. However, computation of $p(\b y_h|\b x_h,\b z_{hi})$ is usually infeasible even if the conditional probability $p(\b y_{hi}|\b U_h=\b u,\b V_{hi}=\b v,\b x_h,\b z_h)$ is obtained by recursion (\ref{eq:recursion}). For this reason, we suggest below a pairwise likelihood based approach.
\section{Pairwise likelihood inference}\label{sec.algo}
In order to make inference on the model parameters, we exploit the
following pairwise log-likelihood:
\begin{eqnarray*}
p\ell(\b\th) &=&
\sum_{h=1}^H\sum_{i=1}^{n_h-1}\sum_{j=i+1}^{n_h}
p\ell_{hij}(\b\th),\\
p\ell_{hij}(\b\th)&=& \log[p(\b y_{hi},\b y_{hj}|\b x_h,\b z_{hi},\b z_{hj})],
\end{eqnarray*}
which recalls the pairwise log-likelihood used by \cite{ren:mol:gey:04}.

Note that, when the dimension of each cluster is two ($n_h=2$, $h=1,\ldots,H$), this function
is the exact log-likelihood of the model, since it is based on the
manifest probability of the responses provided by all the possible pairs of subjects in the same cluster.
\subsection{Computation and maximization of the pairwise likelihood}
In order to efficiently compute the probability $p(\b y_{hi},\b
y_{hj}|\b x_h,\b z_{hi},\b z_{hj})$ as a function of the parameters in $\b\th$, we exploit
recursion (\ref{eq:recursion}) already used for the model
illustrated in Section \ref{sec.LM}. In fact, we have that
\[
p(\b y_{hi},\b
y_{hj}) = p(\tilde{\b y}_{hij}^{(1)},\ldots,\tilde{\b
y}_{hij}^{(T)}|\b x_h,\b z_{hi},\b z_{hj}),
\]
where $\tilde{\b y}_{hij}^{(t)}$ is a realization of the vector
$\tilde{\b Y}_{hij}^{(t)}=(Y_{hi}^{(t)},Y_{hj}^{(t)})\tr$. It may
be simply proved that, for $t=1,\ldots,T$, these vectors follow a
bivariate LM model with covariates since they are conditionally independent given
the latent process $\b W_{hij}^{(1)},\ldots,\b W_{hij}^{(T)}$,
where $\b W_{hij}^{(t)}=(U_h^{(t)},V_{hi}^{(t)},V_{hj}^{(t)})$, and the corresponding
covariates. In particular, this latent process follows a Markov chain with an augmented space of $k=k_1k_2^2$ states indexed by $\b w=(u,v_1,v_2)$. It is simple to see that the initial
probability of state $\b w$ is
\begin{equation}
\phi_{\bl w}=p(\b W_{hij}^{(1)}=\b w)=\la_u\pi_{v_1}\pi_{v_2},\label{eq:init_new}
\end{equation}
whereas, for $t=2,\ldots,T$, transition probability from state $\bar{\b w}=(\bar{u},\bar{v}_1,\bar{v}_2)$ to $\b w$ is
\begin{equation}
\phi_{\bl w|\bar{\bl w}} =p(\b W_{hij}^{(t)}=\b w|\b W_{hij}^{(t-1)}=\bar{\b w})=
\la_{u|\bar{u}}\pi_{v_1|\bar{v}_{1}}\pi_{v_2|\bar{v}_2}.\label{eq:trans_new}
\end{equation}
Moreover, in the case of discrete or categorical response variables, the model assumptions imply that, given $\b W_{hi}^{(t)}=\b w$, the conditional probability of $\tilde{\b y}_{hij}^{(t)}$ is equal to
\begin{equation}
p(\tilde{\b y}_{hij}^{(t)}|\b W_{hi}^{(t)}=\b w,\b x_h,\b z_{hi},\b z_{hj})=
p(y_{hi}^{(t)}|u,v_1,\b x_h^{(t)},\b z_{hi}^{(t)})
p(y_{hj}^{(t)}|u,v_2,\b x_h^{(t)},\b z_{hj}^{(t)}).\label{eq:cond_new}
\end{equation}
A similar expression holds for continuous response variables, based on the corresponding density functions.

In order to compute $p({\b y}_{hi},{\b y}_{hj}|\b x_h,\b z_{hi},\b z_{hj})$, recursion
(\ref{eq:recursion}) is applied with $\b m_i$ substituted by the vector $\tilde{\b m}_{hij}$ having elements $p(\tilde{\b y}_{hij}^{(t)}|\b W_{hi}^{(t)}=\b w,\b x_h,\b z_{hi},\b z_{hj})$ for all $\b w$. Similarly, $\b\pi$ must be substituted by the initial
probability vector $\b\phi$ with elements $\phi_{\bl w}$ and ${\b \Pi}$ by the transition matrix $\b\Phi$ with elements $\phi_{\bl w|\bar{\bl w}}$.

The pairwise log-likelihood $p\ell(\b\th)$ can be maximized by an
EM algorithm having a structure that closely recalls that outlined in Section \ref{sec:EM}. In this case, in particular, the {\em complete data pairwise log-likelihood} is
$$
p\ell^*(\b\th)=\sum_{h=1}^H\sum_{i=1}^{n_h-1}\sum_{j=i+1}^{n_h}
p\ell^*_{hij}(\b\th),
$$
where
\begin{eqnarray}
p\ell^*_{hij}(\b\th)&=& \sum_{\bl w}d_{hij}^{(1)}(\b w)\log(\phi_{\bl w})\nonumber\\
&& + \sum_{t>1}\sum_{\bar{\bl w}}\sum_{\bl w}d_{hij}^{(t)}(\bar{\b w},{\b w})
\log(\phi_{\bl w|\bar{\bl w}})\nonumber\\
&&+\sum_t\sum_{\bl w} d_{hij}^{(t)}(\b w)\log[p(\tilde{\b y}_{hij}^{(t)}|\b W_{hi}^{(t)}=\b w,\b x_h,\b z_{hi},\b z_{hj})].\label{eq:comp_lk}
\end{eqnarray}
In the above expression, $d_{hij}^{(t)}(\b w)$ is a dummy
variable equal to 1 if, at occasion $t$, cluster $h$ is in latent
state $u$, subject $hi$ is in latent state $v_1$, and subject $hj$
is in latent state $v_2$; moreover, we have $d_{hij}^{(t)}(\bar{\b w},{\b w})=d_{hij}^{(t-1)}(\bar{\b w})d_{hij}^{(t)}(\b w)$.

The complete data pairwise log-likelihood may be simply expressed in
terms of the parameters of the proposed multilevel model by
substituting (\ref{eq:init_new}), (\ref{eq:trans_new}), and
(\ref{eq:cond_new}) in the above expression. For instance, the first
component becomes the sum over $u$ of
\begin{equation}
\tilde{d}_{hij}^{(1)}(u)\log[\la_h(u)]+\sum_{v_1}\tilde{d}_{hij}^{(11)}(u,v_1)\log[\pi_{hi}(v_1|u)]+
\sum_{v_2}\tilde{d}_{hij}^{(12)}(u,v_2)\log[\pi_{hj}(v_2|u)],\label{eq:simplif}
\end{equation}
where the variables $\tilde{d}_{hij}^{(1)}(u)$,
$\tilde{d}_{hij}^{(11)}(u,v_1)$, and $\tilde{d}_{hij}^{(12)}(u,v_2)$
are obtained by summing $d_{hij}^{(1)}({\b w})$ over suitable
configurations of ${\b w}$. In a similar way we can express the
other two components involving the transition and the conditional
response probabilities (or densities).

At the E-step of the EM algorithm, the conditional expected value of each dummy variable $d_{hij}^{(t)}({\b w})$ and $d_{hij}^{(t)}(\bar{\b w},{\b w})$ is computed by using the same recursions exploited in the algorithm of \cite{baum:et:al:70}. At the M-step, the model parameters are updated by maximizing the function resulting by substituting the expected values in (\ref{eq:comp_lk}) and exploiting the simplification (\ref{eq:simplif}) and similar simplifications. In any case, the final algorithm is implemented along the same lines as the algorithm implemented by \cite{bart:farc:09}. We make our \verb"R" implementation available to the reader upon request.
\subsection{Model selection and hypothesis testing}
As in \cite{ren:mol:gey:04}, we estimate the variance-covariance matrix of the pairwise likelihood estimator $\hat{\b \th}$, and then obtain standard errors, by the following sandwich formula
\[
\hat{\b\Sigma}(\hat{\b\th})=\hat{\b J}^{-1}\hat{\b K} \hat{\b J}^{-1},
\]
where
$$
\hat{\b J} = -\sum_h \frac{\pa^2 p\ell_h(\hat{\b\th})}{\pa\b\th\pa\b\th\tr},\quad
\hat{\b K} = \sum_h\frac{\pa
p\ell_h(\hat{\b\th})}{\pa\b\th}\frac{\pa p\ell_h(\hat{\b
\th})}{\pa\b\th\tr}, \quad  p\ell_h(\b \th)= \sum_{i=1}^{n_h-1}
\sum_{j=i+1}^{n_h} p\ell_{hij}(\b \th).
$$
We obtain the first derivative of $p\ell_h(\b\th)$ as a by-product
of the EM algorithm. The second derivative, instead, is
obtained by a numerical method.

General results on the asymptotic properties of  the pairwise
likelihood estimator $\hat{\b \th}$ can be derived along the
lines of classical maximum likelihood estimators. However, the
former is expected to be less efficient since it relies on a
restricted amount of information \citep{ren:mol:gey:04}.

In order to deal with model selection, \citet{var:vid:05} suggested CLIC. According to this criterion, the model to be selected is the one which maximizes the following index
\begin{equation}
CLIC = p\ell(\hat{\b\th}) - {\rm tr}(\hat{\b K}\hat{\b J}^{-1}).\label{eq:clic}
\end{equation}
We use this criterion to select the number of states $k_1$ and
$k_2$ of each latent process ${\b U}_h$ at cluster level and ${\b V}_{hi}$
at unit level. Moreover, it can be also used for selecting one
of the possible parametrizations illustrated in Section
\ref{sec.ext}.
\section{Application}\label{sec.appl}
We illustrate the proposed approach by an application based on a dataset on
individual work histories derived from the administrative archives
of the Italian National Institute of Social Security (INPS). We consider a sample of 1,876
employees (both blue-collars and white-collars) from 249 private
Italian firms with 1,000 to 10,000 workers. The subjects,
continuously working in the same firm and aged between 18 and 60 in 1994, were followed for 6 years, from 1994 to 1999. See \citet{bart:nigr:07} for further details.

As already mentioned in Section \ref{se.intro}, the binary response variable of interest is \textit{illness} (equal to 1 if the employee received illness benefits in a certain year and to 0 otherwise). We also consider a set of unit- and cluster-level covariates: \textit{gender} (dummy equal to 1 for woman), \textit{age} in 1994, \textit{area} (Noth-West, North-East, Center, South, or Islands), \textit{skill} (dummy equal to 1 for a blue-collar), \textit{income} (total annual
compensation in thousands of Euros), and \textit{part-time} (dummy
equal to 1 for a part-time employee). Among the covariates we also
include the lagged response.

To this dataset, we fitted the model described in Section \ref{sec.ext} under the constraint that the transition matrices for both processes are tridiagonal with constant off-diagonal elements; see equation (\ref{eq:tri}). We also assume a logistic regression model as in (\ref{eq:cov_psi}) for the conditional probabilities. Then, the unit-level latent
process is expected to capture the propensity (which is not explained by the observed covariates) to get ill of every subject, whereas the cluster-level latent process explains the effect of different firms on the propensity to require illness benefits.

The first step of the analysis is the choice of the number of states for the cluster- and unit-level latent processes, denoted by $k_1$ and $k_2$ respectively. This choice is based on CLIC, which is based on the index defined in (\ref{eq:clic}). The value of this index is reported in Table \ref{tab:clic} for different values of $k_1$ and $k_2$. According to these results we select the model with $k_1=3$ states at cluster level and $k_2=2$ at unit level.

\setlength{\tabcolsep}{5pt}
\begin{table}[ht]
\begin{center}
\caption{\em Values of CLIC for different values of $k_1$ and $k_2$ (in
boldface the largest CLIC value).}\label{tab:clic}
\begin{tabular}{crrrrrr}\hline
  && \multicolumn5c{$k_2$}\\\cline{3-7}
$k_1$ && \multicolumn1c{1} && \multicolumn1c{2} &&
\multicolumn1c{3}\\\hline
1 &&      -30724   &&    -30300    &&   -29972\\
2 &&      -30144   &&    -29773    &&   -29779\\
3 &&      -30018   &&    {\bf -29705}    &&   -29756\\
4 &&      -30001   &&    -29727    &&   -29747\\\hline
\end{tabular}
\end{center}
\end{table}

Table \ref{tab.cov} collects the estimates of the regression parameters obtained with the selected number of states. We note that the probability of receiving illness benefits is positively related to being a blue-collar and to the lagged response, whereas it is negatively related to income and to having a part-time job. The effects of gender, age and age squared are not significant.

\setlength{\tabcolsep}{5pt}
\begin{table}
\caption{\em Estimates of the logistic regression parameters (collected in the vectors $\b\ga$ and $\b\de$) affecting the conditional probabilities.}\label{tab.cov}
\begin{center}
\begin{tabular}{lcrrrrrrrrr}\hline
parameter        && estimate         &&   \multicolumn1c{s.e.} &&
\multicolumn1c{$t$-stat} && \multicolumn1c{$p$-value}\\\hline
intercept        && -3.474 && 1.364 && -2.547 && 0.011\\
gender           &&  0.161 && 0.184 &&  0.876 && 0.382\\
age              && -0.003 && 0.045 && -0.067 && 0.947\\
age$^2/100$      &&  0.038 && 0.060 &&  0.633 && 0.527\\
area: North-East &&  0.145 && 0.257 &&  0.564 && 0.573\\
area: Center     && -0.096 && 0.284 && -0.338 && 0.735\\
area: South      && -0.427 && 0.355 && -1.203 && 0.229\\
area: Islands    && -1.046 && 0.485 && -2.157 && 0.031\\
skill            &&  2.037 && 0.423 &&  4.816 && 0.000\\
income           && -0.200 && 0.035 && -5.714 && 0.000\\
part-time        && -0.795 && 0.338 && -2.352 && 0.019\\
lagged-response  &&  0.600 && 0.172 &&  3.480 && 0.000\\\hline
\end{tabular}

\end{center}
\end{table}

About the distribution of each cluster and unit-level latent process, the estimates of
the initial and transition probabilities are reported in Table
\ref{tab.cluster_proc} and \ref{tab.indiv_proc}. For both processes,
we observe that the states are well separated and the second state
is the one with the highest initial probability. Moreover, the
estimates of the transition matrices show that the cluster-level
latent process has a lower persistence than the unit-level latent process.
\setlength{\tabcolsep}{5pt}
\begin{table}
\begin{center}
\caption{\em Support points and initial and transition probabilities of
each cluster-level latent process.}\label{tab.cluster_proc}
\begin{tabular}{ccccccccccc}\hline
\multicolumn1c{latent} && \multicolumn1c{support} && \multicolumn1c{initial}     && \multicolumn5c{transition}\\
\multicolumn1c{state ($u$)}  && \multicolumn1c{point ($\al_u$)}   &&
\multicolumn1c{probability ($\la_u$)} &&
\multicolumn5c{probabilities ($\la_{u|\bar{u}}$)}\\\hline
1 && 0.000 &&  0.2221 && 0.9130 && 0.0870 && 0.0000\\
2 && 0.444 &&  0.7181 && 0.0870 && 0.8260 && 0.0870\\
3 && 2.931 &&  0.0598 && 0.0000 && 0.0870 && 0.9130\\\hline
\end{tabular}

\end{center}
\end{table}
\setlength{\tabcolsep}{5pt}
\begin{table}
\begin{center}
\caption{\em Support points and initial and transition probabilities of
each unit-level latent process.}\label{tab.indiv_proc}
\begin{tabular}{cccccccccc}\hline
\multicolumn1c{latent} && \multicolumn1c{support} && \multicolumn1c{initial}     && \multicolumn3c{transition}\\
\multicolumn1c{state ($v$)}  && \multicolumn1c{point ($\be_v$)}   &&
\multicolumn1c{probability ($\pi_v$)} &&
\multicolumn3c{probabilities ($\pi_{v|\bar{v}}$)}\\\hline
1 && 0.000 &&  0.4122 && 0.9729 && 0.0271\\
2 && 2.718 &&  0.5878 && 0.0271 && 0.9729\\\hline
\end{tabular}

\end{center}
\end{table}

Finally, we tried to simplify the model selected above by
restricting the transition matrix of each latent process to be
diagonal, so that transition between latent states is not allowed.
In particular, the model in which the transition matrix at
cluster-level is diagonal has a slightly lower value of CLIC
equal to -29,706. On the other hand, the restriction that the
transition matrix at unit-level is diagonal leads to a
strong decrease of CLIC, which is equal to -29,757. We then retain
the model in which latent transition is allowed both at cluster
and unit levels.
\section{Conclusions}
With reference to multilevel longitudinal data, where sample units are collected in clusters, in this paper we propose an approach to account for the unobserved heterogeneity between sample units and between clusters in a dynamic fashion. The approach is based on associating a hidden (or latent) Markov chain to every sample unit and to every cluster. These Markov chains are assumed to be homogeneous and of the first-order, with transition probabilities that may be subjected to suitable constraints. The approach then extends the one proposed by \cite{bart:farc:09}, who proposed a latent Markov model with covariates for longitudinal data (not having a multilevel structure).

The complexity of the model formulated on the basis of the proposed approach does not allow us to make exact likelihood inference on its parameters. Therefore, we adopt a composite likelihood framework for making inference, which is based on considering all the possible pairs of units in every cluster, as suggested by \cite{ren:mol:gey:04} in a simpler context. Within this framework, we also deal with model selection, based on the composite likelihood information criterion \citep{var:vid:05}, and hypothesis testing. In an application based on data about a sample of Italian workers who are employed in different firms, we observed that this composite likelihood approach gives sensible estimates. In this application the response variable is binary, but the approach is completely general in terms of type of response variable, which may be also continuous, discrete, or ordinal.

Possible further developments of the proposed approach may concern the implementation of faster algorithms for the maximization of the pairwise likelihood that we use. In fact, we maximize this function by an Expectation-Maximization (EM) algorithm which is implemented along the same lines as in \cite{bart:farc:09}. However, we think that this maximization may be made much faster by using, after a certain number of EM iterations, a Newton-Raphson algorithm. The implementation of this algorithm is made possible by the availability of the score and the observed information matrix (for the pairwise likelihood function), that we already are able to compute within the present approach.

Finally, another point that deserves attention is the use of alternative forms of composite likelihood for parameter estimation. In particular, in the current form, the adopted pairwise likelihood gives more weight to the data referred to the units belonging clusters having a higher dimension. Then, as suggested by \citet{ren:mol:gey:04}, a weighted version of the pairwise log-likelihood may be more suitable when the clusters are strongly different in terms of dimension. Note, however, that in our application the clusters are not very different in terms of dimension and so, at least in the present case, we do not expect to obtain very different results on the basis of a weighted composite likelihood function.

\bibliography{reference-LM_LM}
\bibliographystyle{apalike}

\end{document}